\documentclass[11pt]{article}

\usepackage{amsmath}
\usepackage{amsfonts}
\usepackage{amssymb}
\usepackage{amsthm}
\usepackage{enumitem}
\usepackage{mathtools}
\usepackage{setspace}
\usepackage{etoolbox}

\newtheorem{theorem}{Theorem}[section]

\newtheorem{lemma}[theorem]{Lemma}
\newtheorem{proposition}[theorem]{Proposition}
\newtheorem{corollary}[theorem]{Corollary}

\renewenvironment{proof}[1][Proof.]{\vspace{-16.5pt} \begin{trivlist}
        \item[\hskip \labelsep {\bfseries #1}]}{\qed \end{trivlist}}

\setitemize{noitemsep, topsep=5.5pt, parsep=5.5pt, partopsep=0pt}

\allowdisplaybreaks
\appto\normalsize{
        \abovedisplayskip=5.5pt plus 2pt minus 2pt
        \belowdisplayskip=5.5pt plus 2pt minus 2pt
        \abovedisplayshortskip=5.5pt plus 2pt minus 2pt
        \belowdisplayshortskip=5.5pt plus 2pt minus 2pt}
\appto\small{
        \abovedisplayskip=5.5pt plus 2pt minus 2pt
        \belowdisplayskip=5.5pt plus 2pt minus 2pt
        \abovedisplayshortskip=5.5pt plus 2pt minus 2pt
        \belowdisplayshortskip=5.5pt plus 2pt minus 2pt}

\newcommand{\gap}{\vspace{11pt}}

\newcommand{\tr}{\operatorname{tr}}

\newcommand{\R}{\mathcal{R}}
\newcommand{\Rn}{\mathcal{R}^n}

\newcommand{\Sn}{\mathcal{S}^n}

\newcommand{\Hn}{\mathcal{H}^n}

\newcommand{\V}{{\cal V}}
\newcommand{\W}{{\cal W}}

\title{\bf Korovkin-type results and doubly stochastic \\
transformations over  
Euclidean Jordan algebras }
\author{
        M. Seetharama Gowda\\
        Department of Mathematics and Statistics\\
        University of Maryland, Baltimore County\\
        Baltimore, Maryland 21250, USA\\
        gowda@umbc.edu
}

\date{\today}

\begin{document}

\maketitle

\begin{abstract}
A well-known theorem of Korovkin asserts that if $\{T_k\}$ is a sequence of positive linear transformations on 
$C[a,b]$ such that $T_k(h)\rightarrow h$ (in the sup-norm on $C[a,b]$) 
for all $h\in \{1,\phi,\phi^2\}$, where $\phi(t)=t$ on $[a,b]$, then $T_k(h)\rightarrow h$ for all $h\in C[a,b]$. In particular,  
if $T$  is a positive linear transformation on $C[a,b]$ such that $T(h)=h$ for all $h\in \{1,\phi,\phi^2\}$, then $T$ is the Identity transformation.
In this paper, we present some analogs of these results over Euclidean Jordan algebras. 
We show that if  $T$ is a positive linear transformation on a Euclidean Jordan algebra $\V$ 
 such that $T(h)=h$ for all $h\in \{e,p,p^2\}$, where $e$ is the unit element in $\V$ and $p$ 
is an  element of $\V$ with  distinct eigenvalues, 
then $T=T^*=I$ (the Identity transformation) on the span of the 
Jordan frame corresponding to the spectral decomposition of $p$; consequently, 
if a positive linear transformation coincides with the Identity 
transformation on a Jordan frame, then it is  doubly stochastic.  
We also present  sequential and weak-majorization versions.

\end{abstract}

\vspace{1cm}
\\

\newpage

%%%%%%%%%%%%%%%%%%%%%%%%%%%%%%%%%%%%%%%%%%%%%%%%%%%%
\section{Introduction}
A well-known theorem of Korovkin \cite{korovkin} asserts that if $\{T_k\}$ is a sequence of positive linear transformations on
the space $C[a,b]$ (of all real-valued continuous functions on the interval $[a,b]$ with sup-norm) such that $T_k(h)\rightarrow h$ 
for all $h\in \{1,\phi,\phi^2\}$, where $\phi(t)=t$ on $[a,b]$, then $T_k(h)\rightarrow h$ for all $h\in C[a,b]$. In particular,
if $T$  is a positive linear transformation on $C[a,b]$ such that $T(h)=h$ for all $h\in \{1,\phi,\phi^2\}$, then $T=I$ (the Identity transformation).
There are numerous generalizations and analogs of Korovkin's theorem in various settings such as Banach function spaces, $C^*$-algebras, etc., see e.g, \cite{altomare-campiti, altomare, priestley}. 
In many of these settings, associativity of the product (as  in $C[a,b]$ and the space of $n\times n$ complex matrices) and an inequality of the form 
$T(h)^2\leq T(h^2)$ (known as Kadison's inequality) are crucially used \cite{uchiyama, priestley}. 
In this paper, we focus on Euclidean Jordan algebras, where, generally, associativity is not available and a Kadison-type inequality is not (yet) known. Here, 
we formulate several Korovkin-type results and make an interesting connection to doubly stochastic transformations. 
\\

Let $\V$ be a Euclidean Jordan algebra with unit element $e$ (see Section 2 for definitions and examples) and   $T$ be a positive linear transformation on it (so $T$ keeps the symmetric cone of $\V$ invariant).  
Let $p$ be an element of $\V$  with  distinct eigenvalues. Our main result,  Theorem \ref{eja korovkin}, asserts that the condition $T(h)=h$ for all $h\in \{e,p,p^2\}$ is equivalent to $T=I$ and also to $T^*=I$ on the span of the Jordan frame corresponding to $p$, where $I$ denotes the Identity transformation on $\V$ and $T^*$ denotes the adjoint of $T$. An immediate consequence is that if a  positive linear transformation coincides with the Identity transformation on a Jordan frame, then it is doubly stochastic (i.e., it is  positive, unital, and trace-preserving).  
The sequential version -- proved as a consequence of our main result -- is as follows: Let $p$ be as above and suppose $\{T_k\}$ is a sequence of positive linear transformations on $\V$ such that $T_k(h)\rightarrow h$ for all $h\in \{e,p,p^2\}$. Then $T_k(h)\rightarrow h$  and $T_k^*(h)\rightarrow h$ for all $h$ in the span of the Jordan frame corresponding to $p$.
Along with the above equality and sequential versions, we also discuss (weak) majorization formulations. We show that under certain conditions, a positive linear transformation $T$ satisfying $T(e)\underset{w}{\prec} e$,  $p\underset{w}{\prec} T(p)$, and
$T(p^2)\underset{w}{\prec}p^2$ coincides with an (algebra) automorphism of $\V$ on the Jordan frame of $p$. We also formulate the problem  of characterizing positive linear transformations $T$ for which the conditions $T(e)\prec e$,  $T(p)\prec p$, and
$T(p^2)\prec p^2$ hold.
\\

An outline of the paper is as follows. We cover some preliminary material in Section 2.
Section 3 deals with   a Korovkin-type result  for matrices and its  weak-majorization modification. 
In Section 4, we present our main result (Theorem \ref{eja korovkin}), describe its connection to Priestley's generalization of 
Korovkin's theorem (\cite{priestley}, Theorem 1.3), and provide some examples. Section 5 deals with  
sequential and weak-majorization versions on Euclidean Jordan algebras.
 %%%%%%%%%%%%%%%%%%%%%%%%%%%%
\section{Preliminaries}
In $\Rn$, vectors are considered as either column vectors or row vectors depending on the context. For any $x=(x_1,x_2,\ldots, x_n)\in \Rn$, 
we write $x^2:=(x_1^2,x_2^2,\ldots, x_n^2).$
We  write  $e$ for the vector of ones (reserving the same symbol for the  unit element in a general Euclidean Jordan algebra, see below). 
We say that a real $n\times n$ matrix $A$ is {\it nonnegative} if all its entries are nonnegative; it is {\it unital} if $Ae=e$  and {\it subunital} if $Ae\leq e$.
A nonnegative unital matrix is said to be (row) {\it stochastic}. 
A nonnegative matrix $A$ with both $A$ and $A^T$ unital is said to be {\it doubly stochastic.}

\gap

For a vector $x\in \Rn$ with entries/components $x_1,x_2,\ldots, x_n$, let $x^\downarrow $ denote the vector obtained by rearranging the 
of entries of $x$ in a decreasing manner (so that $x_1^\downarrow \geq x_2^\downarrow \geq \cdots\geq x_n^\downarrow$). Clearly  $x^\downarrow =Px$ for some permutation matrix $P$ and $(Ex)^\downarrow =x^\downarrow $ for every
 permutation matrix $E$. (Recall that
a permutation matrix is obtained by permuting the rows/columns of the Identity matrix.)
 It is known (\cite{bhatia}, page 29) that
$$x_1+x_2+\cdots+x_m\leq x_1^\downarrow +x_2^\downarrow+\cdots+x_m^\downarrow$$
for all $m=1,2,\ldots, n$. Given $x,y\in \Rn$, we say that $x$ is weakly-majorized by $y$ and write $x\underset{w}{\prec} y$ if
$x^\downarrow_1+x^\downarrow_2+\cdots+x^\downarrow_k\leq y^\downarrow_1+y^\downarrow_2+\cdots+y^\downarrow_k$ for each $k\in \{1,2,\ldots, n\}$. If, additionally, equality holds for $k=n$, we say that $x$ is majorized by $y$ and write $x\prec y$.  By a well-known theorem of Hardy-Littlewood-Polya (\cite{bhatia}, Theorem II.1.10), $x\prec y$ if and only if $x=Dy$ for some doubly stochastic matrix $D$. Moreover, by Birkhoff's theorem (\cite{bhatia}, Theorem II.2.3), every doubly stochastic matrix is a convex combination of permutation matrices. We note one useful property of majorization: $x\prec y\Rightarrow f(x)\leq f(y)$ for any real-valued convex function $f$ on $\Rn$. 
\gap

The standard material on Euclidean Jordan algebras given below can be found in \cite{faraut-koranyi, gowda-sznajder-tao}. A {\it Euclidean Jordan algebra} is a finite dimensional real inner product space $(\V, \langle\cdot,\cdot\rangle)$ together with a bilinear product (called the Jordan product) $(x,y)\rightarrow x\circ y$ satisfying the following properties:
\begin{itemize}
\item [$\bullet$] $x\circ y=y\circ x$,
\item [$\bullet$] $x\circ (x^2\circ y)=x^2\circ (x\circ y)$, where $x^2=x\circ x$, and
\item [$\bullet$] $\langle x\circ y,z\rangle=\langle x,y\circ z\rangle.$
\end{itemize}
In such an algebra, there is the `unit element' $e$ with the property $x\circ e=x$ for all $x$.
In $\V$,
$$K=\{x\circ x:\,x\in \V\}$$
is called the {\it symmetric cone} of $\V$. It is a self-dual cone.

\gap

The space $\Rn$ is a Euclidean Jordan algebra  under the componentwise product and the usual inner product. In this algebra, the symmetric cone is the nonnegative orthant. 
Any (nonzero) Euclidean Jordan algebra is a direct product/sum
of simple Euclidean Jordan algebras and every simple Euclidean Jordan algebra is isomorphic to one of five algebras,
three of which are the algebras of $n\times n$ real/complex/quaternion Hermitian matrices. The other two are: the algebra ${\cal O}^3$ of $3\times 3$ octonion Hermitian matrices and the Jordan spin algebra ${\cal L}^n$.
In the algebras $\Sn$ (of all $n\times n$ real symmetric matrices) and $\Hn$ (of all $n\times n$ complex Hermitian matrices), the Jordan product and the inner product are given, respectively, by
$$X\circ Y:=\frac{XY+YX}{2}\quad\mbox{and}\quad \langle X,Y\rangle:=\tr(XY),$$
where the trace of a real/complex matrix is the sum of its diagonal entries.

\gap

Let $\V$ be a Euclidean Jordan algebra. A nonzero element $c$ in $\V$ is an {\it idempotent} if $c^2=c$; it is a {\it primitive idempotent} if it is not the sum of two other idempotents. {\it
A Jordan frame $\{e_1,e_2,\ldots, e_n\}$ in $\V$ consists of
primitive idempotents that are mutually orthogonal  (equivalently, $e_i\circ e_j=0$ when $i\neq j$) with sum equal to the unit element.} All Jordan frames in $\V$ have the same number of elements, called the rank of $\V$.
{\it Let the rank of $\V$ be $n$.} According to the {\it spectral decomposition
theorem} \cite{faraut-koranyi}, any element $x\in \V$ has a decomposition
\begin{equation} \label{spectral decomposition}
x=x_1e_1+x_2e_2+\cdots+x_ne_n,
\end{equation}
where the real numbers $x_1,x_2,\ldots, x_n$ are (called) the eigenvalues of $x$ and
$\{e_1,e_2,\ldots, e_n\}$ is a Jordan frame in $\V$. (An element may have decompositions coming from different Jordan frames, but the eigenvalues remain the same. However, if all the eigenvalues are distinct, then, up to permutation, there is only one spectral decomposition, see \cite{faraut-koranyi}, Theorem III.1.1.)
For notational simplicity, we write the above spectral decomposition (\ref{spectral decomposition}) in the form $x=r*{\cal E}$, where $r=(x_1,x_2,\ldots,x_n)$ and ${\cal E}:=\{e_1,e_2,\ldots, e_n\}$.

\gap

{\it For any $x\in \V$, let $\lambda(x)$ denote the vector of eigenvalues of $x$ written in the decreasing order.}
Then, we can always write the spectral decomposition of any $x\in \V$ in the form  $x=\lambda_1(x)f_1+\lambda_2(x)f_2+\cdots+\lambda_n(x)f_n=\lambda(x)*{\cal F}$  relative to a Jordan frame ${\cal F}=\{f_1,f_2,\ldots, f_n\}$. 

\gap

Given $a\in \V$, we define linear transformations $L_a$ and $P_a$ (called the quadratic representation of $a$) on $\V$ 
by
$$L_a(x):=a\circ x\quad\mbox{and}\quad P_a(x):=2a\circ (a\circ x)-a^2\circ x\,\,(x\in \V).$$
We say that elements $a,b\in \V$ {\it operator commute} if the transformations $L_a$ and $L_b$ commute. It is known, see \cite{faraut-koranyi}, Lemma X.2.2, that 
{\it $a$ and $b$ operator commute if and only if $a$ and $b$ have their spectral decompositions with respect to the same Jordan frame.
}

\gap

For any $x\in \V$ with eigenvalues $x_1,x_2,\ldots, x_n$, the {\it trace} of $x$ is defined by
$$\tr(x):=x_1+x_2+\cdots+x_n.$$

It is  known that $(x,y)\mapsto \tr(x\circ y)$ defines another inner product on $\V$ that is compatible with the Jordan product.
We let 
$$\langle x,y\rangle_{tr}:=\tr(x\circ y)$$
and call this, the {\it trace inner product.}
When we replace the given inner product by the trace inner product, Jordan frames as well as the eigenvalues of an element remain the same. 
One advantage is: In the trace  inner product, the norm of any primitive idempotent is one and so any Jordan frame
in $\V$ is an orthonormal set. Additionally,
$\tr(x)=\langle x,e\rangle_{tr} \quad\mbox{for all}\,x\in \V.$

\gap

We use the notation $x\geq 0$ ($x> 0$) when $x\in K$ (respectively, interior of $K$) or, equivalently, all the eigenvalues of $x$ are nonnegative (respectively, positive); when $x>0$, we say that $x$ is a positive element. 
We also write $x\leq y$ in $\V$ when $y-x\geq 0$. 
{\it Since $K$ is self-dual, we see that $x\geq 0$ if and only if $\langle x,y\rangle \geq 0$ for all $y\geq 0$.
} For any $x\in \V$ with spectral decomposition $x=x_1e_1+x_2e_2+\cdots+x_ne_n,$,
we define $x^+:=x_1^+e_1+x_2^+e_2+\cdots+x_n^+e_n,$ and $x^-:=x^+-x$ so that $x^+,x^-\in K$ and $x=x^+-x^-$. (Here, for any real number $\lambda$, $\lambda^+:=\max\{\lambda,0\}$.)
\\

We record one useful  consequence of the well-known Hirzebruch's min-max theorem \cite{hirzebruch}: 
\begin{equation}\label{hirzebruch}
x\leq y\Rightarrow \lambda(x)\leq \lambda(y).
\end{equation}

\gap

Given a Jordan frame  $\{e_1,e_2,\ldots, e_n\}$,
we have the {\it Peirce orthogonal decomposition} (\cite{faraut-koranyi}, Theorem IV.2.1):
$\V=\sum_{i\leq j}\V_{ij},$
where $\V_{ii}:=\{x\in \V: x\circ e_i=x\}=\R\,e_i$ and for $i<j$, $\V_{ij}:=\{x\in \V: x\circ e_i=\frac{1}{2}x=x\circ e_j\}.$
Then, for any $x\in \V$, we have
\begin{equation} \label{long peirce decomposition}
x=\sum_{i\leq j}x_{ij}=\sum_{i=1}^{n}x_ie_i+\sum_{i<j}x_{ij}\quad\mbox{with}\quad  x_i\in \R\,\,\mbox{and}\,\,x_{ij}\in \V_{ij}.
\end{equation}

Let $\{e_1,e_2,\ldots, e_n\}$ be a (fixed) Jordan frame in $\V$. For arbitrary $x,y\in \V$, consider the corresponding Peirce decompositions
$x=\sum_{i=1}^{n}x_ie_i+\sum_{i<j}x_{ij}$ and $y=\sum_{i=1}^{n}y_ie_i+\sum_{i<j}y_{ij}$. Define the real symmetric matrix
$$x\Delta y:=\sum_{i=1}^{n} x_iy_i||e_i||^2E_{ii}+\frac{1}{2}\sum_{i<j} \langle x_{ij},y_{ij}\rangle E_{ij},$$
where $E_{ij}$ is the $n\times n$ matrix with $1$s in the $(i,j)$ and $(j,i)$ slots and zeros elsewhere.
It has been proved in \cite{gowda-tao}, Theorem 8, that $x\Delta y$ is a (symmetric) positive semidefinite matrix when $x, y\geq 0$. Since
$$x\Delta x=\sum_{i=1}^{n} x_i^2||e_i||^2 E_{ii}+\frac{1}{2}\sum_{i<j} ||x_{ij}||^2 E_{ij},$$
and all (in particular, $2 \times 2$ and $1\times 1$) principal minors of a real symmetric positive semidefinite matrix are nonnegative, we see that for $x\geq 0$ and $i<j$,
$$||x_{ij}||^2\leq 2 x_ix_j||e_i||\,||e_j||.$$

(Note: When $\V$ carries the trace inner product, this inequality reduces to $||x_{ij}||^2\leq 2 x_ix_j$, see \cite{faraut-koranyi}, Page 80.) We record a useful consequence:

\begin{proposition}\label{zero diagonal entry}$\,\,$
{\it Suppose $x\geq 0$ in $\V$ and let $x=\sum_{i=1}^{n}x_ie_i+\sum_{i<j}x_{ij}$ be its Peirce decomposition relative to a  given
Jordan frame $\{e_1,e_2,\ldots, e_n\}$. If $x_i=0$ for some $i$, 
then $x_{il}=0$ for $l>i$ and $x_{li}=0$ for $l<i$.
}
\end{proposition}

\gap

{\it A linear transformation $T:\V\rightarrow \V$ is said to be positive if
$x\geq 0\Rightarrow T(x)\geq 0$
and unital if $T(e)=e$.
$T$ is said to be {\it trace-preserving} if $tr(T(x))=tr(x)$ for all $x\in \V$. A positive unital trace-preserving transformation is said to be doubly stochastic.
}

Note that positivity means that $T(K)\subseteq K$, where $K$ is the symmetric cone of $\V$. By the self-duality of $K$, if $T$ is  positive, then so is the adjoint $T^*$ (defined by the condition $\langle T(x),y\rangle=\langle x,T^*(y)\rangle$ for all $x,y\in \V$). Writing $T_{tr}^*$ for the adjoint of $T$ relative to the trace inner product, we see that  $T$ is trace-preserving if and only if $T_{tr}^*(e)=e$.
\\

{\it A linear transformation $\phi:\V\rightarrow \V$ is an (algebra) automorphism if it is  bijective and $\phi(x\circ y)=\phi(x)\circ \phi(y)$ for all $x,y\in \V$.} An automorphism $\phi$ maps a Jordan frame to a Jordan frame and so $\lambda(\phi(x))=\lambda(x)$ for all $x$. Moreover, 
for any $p\in \V$, $\phi(p)$ and $\phi(p^2)$ have their spectral decompositions with respect to the same Jordan frame.
We observe that automorphisms are doubly stochastic. We mention two  more
 useful results:
On a simple Euclidean Jordan algebra, any Jordan frame can be mapped onto any another by an automorphism, see \cite{faraut-koranyi}, Theorem IV.2.5. Also, 
every 
automorphism on $\Hn$ is of the form 
$\phi:X\mapsto UXU^*$ for some unitary matrix $U$.

\gap

{\it We define majorization in $\V$  by: $x\prec y$ in $\V$ if  $\lambda(x)\prec \lambda(y)$  in $\Rn$. Likewise, $x\underset{w}{\prec} y$ if 
$\lambda(x)\underset{w}{\prec} \lambda(y)$  in $\Rn$.}

\gap

We have the following result from \cite{gowda-doubly stochastic}:

\begin{theorem}\label{ds theorem from gowda}$\,\,$
{\it 
For $x, y \in \V$, consider the following statements:
\begin{itemize}
\item [$(a)$] $x = T(y)$, where $T$ is a convex combination of  automorphisms of $\V$.
\item [$(b)$] $x=T(y)$, where $T$ is doubly stochastic on $\V$.
\item [$(c)$] $x\prec y$ in $\V$.
\end{itemize}
Then, $(a)\Rightarrow (b)\Rightarrow (c)$. Furthermore, reverse implications hold when $\V$ is $\Rn$ or simple.
}
\end{theorem}

\gap{}

Along with the above, we  mention a
result of Jeong and Gowda (\cite{jeong-gowda}, Lemma 2):
{\it  $T$  is doubly stochastic
if and only if $T(x)\prec x$ for all $x\in \V$.} For some results related to weak-majorization, we refer to \cite{jeong-jung-lim}.

\gap

Throughout this paper, depending on the context, $I$ denotes either the Identity matrix  or
the Identity transformation (on a vector space).

%%%%%%%%%%%%%%%%%%%%%%%%%%%%%%%%%%%%%%%%%%%%%%%%%%%%%%%%
\section{Results over  $\Rn$}
Our first result is   stated in the setting of  the algebra $\Rn$.
While it can be derived from  
known results such as  Theorem 1.3 in \cite{priestley}, 
for completeness, we provide a simple and direct proof. 

\gap{}

\begin{theorem}\label{matrix korovkin}$\,\,$
{\it 
Suppose  $A\in \R^{n\times n}$ is a nonnegative matrix 
such that $Ah=h$ for all $h\in \{e,p,p^2\}$, where $p\in \Rn$  is a vector with distinct entries. Then $A=I$.
}
\end{theorem}

\gap{}

\begin{proof} Let $A=[a_{ij}]$ and $p=(p_1,p_2,\ldots, p_n)$. For any fixed $i$, we  claim that the $i$th  row of $A$, namely, $(a_{i1},a_{i2},\ldots,a_{in})$ has $1$ in the $i$th slot and zeros elsewhere. 
We observe that the entries of this row are nonnegative.  
The given conditions $Ae=e$, $Ap=p$, and $Ap^2=p^2$ imply that 
$(Ae)_i=1$, $(Ap)_i=p_i$, and $(Ap^2)_i=p_i^2$. Then, from the convexity of 
the function $t\mapsto t^2$ on $\R$, 
$$p_i^2=\Big (\sum_{k=1}^{n} a_{ik}p_k\Big )^2\leq \sum_{k=1}^{n} a_{ik}p_k^2=p_i^2.$$
Consequently, since the entries of $p$ are distinct, by the strict convexity of the function $t\mapsto t^2$, only one $a_{ik}$ can be nonzero. From $\sum_{k=1}^{n}a_{ik}=1$ and $\sum_{k=1}^{n} a_{ik}p_k=p_i$, we see that $k=i$ and $a_{ii}=1$. This proves our claim. Thus, $A=I$. 
 
\end{proof}

\gap
In our next result, we replace the equality $Ah=h$ by an appropriate weak-majorization inequality.

\gap{}

\begin{theorem}\label{weak majorization korovkin}$\,\,$
{\it
Let  $A\in \R^{n\times n}$ be a nonnegative matrix and $p\in \Rn$  be a positive vector with distinct entries.
If $Ae\underset{w}{\prec} e$, $p\underset{w}{\prec} Ap$, and $Ap^2\underset{w}{\prec} p^2$, then  $A$ is a permutation matrix; additionally, if the entries of $p$ and $Ap$ are decreasing, then $A=I$. 
}
\end{theorem}
\gap

Note: In this result we assume that $p$ has positive entries. Without this assumption, the result may not hold. For example, over $\R^2$, let $p$ be a vector with entries $1$ and $-1$, and $A$ be the matrix with 
rows $(1,0)$ and $(\frac{1}{2}, \frac{1}{4})$.

\gap

\begin{proof}
Suppose $Ae\underset{w}{\prec} e$, $p\underset{w}{\prec} Ap$, and $Ap^2\underset{w}{\prec} p^2$. From $Ae\underset{w}{\prec} e$,  we see that $Ae\leq e$, that is, $A$ is subunital. Let 
$q:=p^\downarrow$ so that the {\it entries of $q$ are positive and strictly decreasing}. Then, $q^2=(p^2)^\downarrow$. Let $p=E_1q$ and $Ap=E_2(Ap)^\downarrow$  
for some permutation matrices $E_1$ and $E_2$. With $B:=E_2^{-1}AE_1$, we verify that 
\begin{center}
{\it $B$ is nonnegative, $Be\leq e$, $q\underset{w}{\prec} Bq$, and $Bq^2\underset{w}{\prec} q^2$.  }
\end{center}
Additionally, $Bq=(Ap)^\downarrow$ so {\it $Bq$ has decreasing entries}.
We now claim  that $B=I$.
\\

Consider the first row of $B$. As $q$ and $Bq$ have decreasing entries, from  $q\underset{w}{\prec} Bq$ we have $q_1\leq (Bq)_1$. As $B$ is nonnegative, $Be\leq e$ and $q_1$ is the largest entry in $q$,   
by the convexity of the  function $t\mapsto t^2$ on $\R$,  
\begin{equation}\label{inequality for the first coordinate}
q_1^2\leq \Big  (\sum_{j=1}^{n}b_{1j}q_j\Big )^2 \leq \sum_{j=1}^{n} b_{1j}q_j^2\leq \sum_{j=1}^{n} b_{1j}q_1^2\leq q_1^2.
\end{equation}
From  the ensuing equality, we have $b_{1j}q_j^2=b_{1j}q_1^2$ for all $j$. As the entries of $q$ are positive and distinct, $b_{1j}=0$ for all $j\neq 1$ and,
 from (\ref{inequality for the first coordinate}), $b_{11}=1$. Thus, the first row of $B$ is 
$(1,0,0,\ldots, 0)$. In particular, $q_1=(Bq)_1$. We use induction to show that the $k$th row of $B$ is of the form $(0,\ldots, 0,1,0,\ldots, 0)$ with $1$ in the $k$th slot. 
Assume that this statement holds for all indices in  $\{1,2,\ldots, k\}$, where $k<n$. (From the above argument, this statement holds for $k=1$.) We show that the statement holds of $k+1$. To simplify the notation, let
$l=k+1$, $x:=Bq^2$ and $y:=q^2$. \\
From the induction hypothesis (by the form of $B$), $q_i=(Bq)_i$ for all $i=1,2,\ldots, k$. Then, from $q\underset{w}{\prec} Bq$ and the fact that the entries of $q$ and $Bq$ are decreasing, we have
$$q_l\leq (Bq)_l.$$ From $x=Bq^2\underset{w}{\prec} q^2=y$, we have, for all $m=1,2\ldots, n$, 
\begin{equation}\label{two inequalities}
x_1+x_2+\cdots+x_m\leq x_1^\downarrow+x_2^\downarrow +\cdots+x_m^\downarrow\leq y_1+y_2+\cdots +y_m,
\end{equation}
where the second inequality is due to the fact that the entries of $y$ are decreasing. From the form of the first $k$ rows of $B$, we have 
$x_i=y_i$ for all $i=1,2\ldots, k$; by successively 
putting $m=1,2,\ldots, k$ in (\ref{two inequalities}), we get 
$x_i=x_i^\downarrow =y_i$ for all $i=1,2\ldots, k$. By putting $m=l\, (=k+1)$ in 
(\ref{two inequalities}), we get
$$x_l\leq x_l^\downarrow\leq y_l,\,\mbox{that is,}\,\,(Bq^2)_l\leq (Bq^2)_l^{\downarrow}\leq q^2_l.$$
Since $q_l\leq (Bq)_l$ with $B$ nonnegative and $Be\leq e$,  
by the convexity  of the function $t\mapsto t^2$ on $\R$, 
$$q_l^2\leq  \Big (\sum_{j=1}^{n}b_{lj}q_j\Big )^2 \leq \sum_{j=1}^{n}b_{lj}q_j^2=(Bq^2)_l\leq (q^2)_l.$$
Then, by the ensuing equality and the  
 strict convexity of function $t\mapsto t^2$, we get 
$b_{lj}=0$ for all $j\neq l$ and $b_{ll}=1$. This proves that our induction statement holds for $l$ $(=k+1)$. We conclude that $B=I$.
Now, $I=B=E_2^{-1}AE_1$ implies that $A$ is a product of permutation matrices, hence a permutation matrix.\\ 
Finally, suppose $p$ and $Ap$ have decreasing entries. Since $A$ is a permutation matrix, it follows that $A$ must be the Identity matrix. 
(This can also be seen by letting 
$q=p$ and  $E_1=E_2=I$ in the above proof so that $A=B=I$.) This completes the proof. 
\end{proof} 

\gap

We state two immediate consequences.

\begin{corollary} \label{simple corollary}$\,\,$
{\it Suppose $A\in \R^{n\times n}$ is nonnegative and $p$ is a positive vector with distinct entries. If $(Ah)^\downarrow=h^\downarrow$ for all $h\in \{e,p,p^2\}$, then  $A$ is a permutation matrix; additionally if  the entries of $p$ and $Ap$ are decreasing, then $A=I$.
}
\end{corollary}

\begin{corollary}$\,\,$
{\it Suppose $A\in \R^{n\times n}$ is doubly stochastic (so that $Ax\prec x$ for all $x\in \Rn$). If $(Ap)^\downarrow=p^\downarrow$ for some $p\in \Rn$ with positive and distinct entries, 
then $A$ is a permutation matrix; additionally, if the entries of $p$ and $Ap$ are decreasing, then $A=I$. 
}
\end{corollary}

\gap{}

We remark that the latter corollary can also be deduced from Birkhoff's theorem via the strict convexity of the Euclidean norm.

\gap

In reference to the above theorem, one may ask if the condition $p\underset{w}{\prec} Ap$ can be replaced by $Ap\underset{w}{\prec} p$ to make $A$ (at least) doubly stochastic.  
A simple example (such as the $2\times 2$ matrix with rows $(0,1)$ and $(0,1)$ with $p$ having entries $1$ and $2$) shows that 
 for a nonnegative matrix, the conditions $Ae\underset{w}{\prec} e$, $Ap\underset{w}{\prec} p$, and $Ap^2\underset{w}{\prec} p^2$, need not imply that $A$ 
a doubly stochastic.  What if we  replace weak-majorization inequalities by majorization ones? As the answer is  unclear, we pose the following.

\gap

\noindent{\bf Problem:} {\it Let $p\in \Rn$ be a positive vector with distinct entries. Consider the compact convex set
$$\Omega_p:=\{A\in \R^{n\times n}: \,\,A\,\,\mbox{is nonnegative}\,\,\mbox{and}\,\, Ah\prec h\,\,\mbox{for all}\,\,h\in \{e,p,p^2\}\}.$$
Is every matrix in this set doubly stochastic? If not, what  are the extreme points of this set? 
}\\

When $Ah\prec h$ for some $h$, we have $\langle Ah,e\rangle =\langle h,e\rangle$, that is, $\langle A^Te-e,h\rangle =0$. Hence, if $A$ is nonnegative and
$p$ is a positive vector with distinct entries, then the condition $Ah\prec h$ for all $h$ in (the basis) $\{e,p,p^2,\ldots, p^{n-1}\}$ implies that $A^Te=e$, that is, $A$ is doubly stochastic.
In particular, if $A\in \R^{n\times n}$ with $n\leq 3$, the above problem has an affirmative answer, that is, every $A$ in $\Omega_p$ is doubly stochastic.
The answer for $n\geq 4$ is unclear.

\gap

%%%%%%%%%%%%%%%%%%%%%%%%%%%%%%%%%%%%%%%%%%%%%%%%%%%%%%%%%%%%%%%
\section{Equality versions over general Euclidean Jordan algebras}
 Throughout this section, we assume that 
$\V$ is a Euclidean Jordan algebra of rank $n$ with unit element $e$.

\gap

Before proving our  equality/identity version of Korovkin's theorem over  Euclidean Jordan algebras, we provide a simple example to show
 that the direct analog of Theorem \ref{matrix korovkin} is false.

\gap

\noindent{\bf Example 1} Let $\V=\Hn$ and $T$ be the transformation that takes a matrix  $X\in \V$ to the corresponding diagonal matrix, that is,
$$T(X):=\mbox{Diag} (X),$$
where $\mbox{Diag}(X)$ is the diagonal matrix whose diagonal is that of $X$.
Then, $T$ is  linear, positive, and unital. Moreover, $T(H)=H$ for all $H\in \{I,P,P^2\}$, where $I$ is the Identity matrix, $P$ is any diagonal matrix with  distinct diagonal entries. Yet,  $T$ does not coincide with the Identity transformation  on $\Hn$.  
It is interesting to observe that $T$ is trace-preserving and, hence,  doubly stochastic. It follows, for example, from Theorem \ref{ds theorem from gowda}, that the diagonal of a Hermitian matrix is majorized by the eigenvalue vector of that matrix -- this is the well-known Schur's theorem in matrix theory. 

\gap

In preparation for the main theorem, we present a lemma. Here, we let $\delta_{ij}$  denote Kronecker's delta function. 

\begin{lemma}\label{basic lemma}$\,\,$
{\it Suppose $T:\V\rightarrow \V$  is a positive linear transformation and $\{e_1,e_2,\ldots, e_n\}$ is a Jordan frame such that
\begin{equation}\label{delta equation}
\langle T(e_j),e_i\rangle=||e_i||^2\,\delta_{ij}\quad (1\leq i,j\leq n).
\end{equation}
Then, $T(e_k)=T^*(e_k)=e_k$ for all $k$; moreover,  $T$ is doubly stochastic.
}
\end{lemma}

\gap{}

\begin{proof}
Let $I$ denote the Identity transformation on $\V$. Fix any $k\in \{1,2,\ldots, n\}$ and  let $b:=T(e_k)$. Consider the Peirce decomposition of $b$ relative to the Jordan frame $\{e_1,e_2,\ldots, e_n\}$:
$$b=\sum_{i=1}^{n}b_{i}e_i+\sum_{1\leq i<j\leq n}b_{ij}.$$
Now fix any index $i$,  $i\neq k$. By (\ref{delta equation}) and
the orthogonality of the individual terms in the Peirce decomposition,
$$0=\langle T(e_k),e_i\rangle =\langle b,e_i\rangle =b_{i}\,||e_i||^2.$$ 
Since $b\geq 0$ (due to the positivity of $T$), from Proposition \ref{zero diagonal entry}, we must have $b_{il}=0$ for all $l> i$ and
$b_{li}=0$ for all $l<i$.
So, in the Peirce decomposition of $b$, only one term survives.
Hence, $b=b_{k}e_k$. Since $||e_k||^2=\langle T(e_k),e_k\rangle =\langle b,e_k\rangle=b_k||e_k||^2$, we must have $b_k=1$. Thus, $b=e_k$, proving the equality $T(e_k)=e_k$. As $k$ is arbitrary, $T=I$ on the Jordan frame $\{e_1,e_2,\ldots, e_n\}$ and  on the span of $\{e_1,e_2,\ldots, e_n\}$. \\ 
We now show that $T^*=I$ on this span.
 As $T$ is positive, $T(K)\subseteq K$. Since $K$ is self-dual, $T^*(K)\subseteq K$, so $T^*$ is also positive.
Since the condition $\langle T(e_j),e_i\rangle =||e_i||^2\,\delta_{ij}$ is the same as $\langle T^*(e_i),e_j\rangle =||e_j||^2\delta_{ij}$, from the above proof we see that
$T^*(e_k)=e_k$ for all $k$; hence $T^*=I$ on the span of $\{e_1,e_2,\ldots, e_n\}$; in particular, $T^*(e)=e$. 
As $T$ is positive and $T(e)=e$, to show that $T$ is doubly stochastic, we need only show that $T$ is trace-preserving, that is, $T_{tr}^*(e)=e$, where
$T_{tr}^*$ is the adjoint of $T$ relative to the trace inner product. From $T(e_k)=e_k$ for all $k$, we see that 
$$\langle T(e_j),e_i\rangle_{tr}=\langle e_j,e_i\rangle_{tr}=||e_i||_{tr}^2\,\delta_{ij}\quad (1\leq i,j\leq n).$$
By what has been proved earlier (applied to the trace inner product), $T_{tr}^*(e_k)=e_k$ for all $k$ and so, $T_{tr}^*(e)=e$. Thus, $T$ is doubly stochastic.
\end{proof}

\gap

We now state our main theorem.

\begin{theorem}\label{eja korovkin}$\,\,$
{\it 
Suppose $\V$ is a Euclidean Jordan algebra of rank $n$ and $T:\V\rightarrow \V$  is a positive linear transformation. 
Let $\{e_1,e_2,\ldots, e_n\}$ be a Jordan frame and $p=p_1e_2+p_2e_2+\cdots+p_ne_n$, where $p_i$s are  distinct. 
Then, the following statements are equivalent:
\begin{itemize}
\item [$(a)$] $T(h)=h$ for all $h\in \{e,p,p^2\}$.
\item [$(b)$] $T(h)=h$ for all $h\in \mbox{span}\{e_1,e_2,\ldots, e_n\}$.
\item [$(c)$] $T^*(h)=h$ for all $h\in \mbox{span}\{e_1,e_2,\ldots, e_n\}$.
\item [$(d)$] $T^*(h)=h$ for all $h\in \{e,p,p^2\}$.
\end{itemize}
Moreover, under any of the above conditions,
$T$ is doubly stochastic; hence $T(x)\prec x$ for all $x\in \V$.
}
\end{theorem}

\gap{}

\begin{proof} $(a)\Rightarrow (b)$, $(a)\Rightarrow (c)$: Suppose $(a)$ holds. As $p=p_1e_1+p_2e_2+\cdots+p_ne_n$, we have $p^2=p^2_1e_1+p^2_2e_2+\cdots+p^2_ne_n$.
Consider the matrix $A=[a_{ij}]$, where $a_{ij}:=\frac{1}{||e_i||^2}\langle T(e_j),e_i\rangle.$ As $T$ is positive and the symmetric cone $K$ is self-dual, we see that $A$  is a nonnegative matrix.
Since the sum of all $e_j$s is $e$ with $T(e)=e$ and $\langle e,e_i\rangle=||e_i||^2$, we see that 
$A$  is stochastic. From $T(p)=p$ and $T(p^2)=p^2$, we  have, for all $i$, 
$$\sum_{j=1}^{n} a_{ij}p_j=p_i\quad\mbox{and}\quad \sum_{j=1}^{n}a_{ij}p_j^2=p_i^2.$$
Thus, $A$ satisfies the conditions of Theorem \ref{matrix korovkin}. It follows that $A$  is the Identity matrix. Hence,
$\langle T(e_j),e_i\rangle=||e_i||^2\,\delta_{ij},$ for all  $1\leq i,j\leq n$.
From the above lemma, $T(e_k)=T^*(e_k)=e_k$ for all $k$. We see that $T=T^*=I$ 
on the span of $\{e_1,e_2,\ldots,e_n\}$. Thus we  have $(b)$ and $(c)$. \\
$(b)\Rightarrow (a)$, $(c)\Rightarrow (d)$: These are obvious, as $p,p^2\in \mbox{span}\{e_1,e_2,\ldots, e_n\}$.\\
$(d)\Rightarrow (c)$, $(d)\Rightarrow (b)$: Suppose $(d)$ holds. Since $T$ is positive and the symmetric cone $K$ is self-dual, $T^*$ is also positive. Then,
applying the implications $(a)\Rightarrow (b)$, $(a)\Rightarrow (c)$ to $T^*$, we see that $(c)$ and $(b)$ hold. \\
From the above, we see that the stated conditions are all equivalent. Now suppose $(b)$ holds. Then, 
condition (\ref{delta equation}) holds. From the above lemma, $T$ is doubly stochastic.
By Theorem \ref{ds theorem from  gowda}, $T(x)\prec x$ for all $x\in \V$. 
\end{proof}

\gap

\begin{corollary} \label{T equals an automorphism}$\,\,$
{\it 
If a positive linear transformation  coincides with an  automorphism on a Jordan frame, then it is doubly stochastic.
}
\end{corollary}
\gap{}

\begin{proof}
Let $T=\phi$ on a Jordan frame, where $T$ is a positive linear transformation and $\phi$ is an automorphism. Then, using the properties of $\phi$ and $\phi^{-1}$, we see that $\phi^{-1}\circ T$ is a positive linear transformation satisfying the condition $(b)$ of the above theorem. 
Hence $\phi^{-1}\circ T$ is doubly stochastic. As $\phi$ is doubly stochastic, we see that the composition $\phi\circ (\phi^{-1}\circ T)$ is also doubly stochastic. 
Hence, $T$ is doubly stochastic. 
\end{proof}

\gap

\noindent{\bf Remarks.} 
In the setting of a $C^*$-algebra ${\cal A}$, Kadison's inequality  \cite{kadison} asserts that  if a linear transformation $T:{\cal A}\rightarrow {\cal A}$ is positive and $T(I)\leq I$, then
$$T(X)^2\leq T(X^2)$$
for all self-adjoint elements $X$ in ${\cal A}$. (Uchiyama's elementary proof of Korovkin's theorem \cite{uchiyama} uses this inequality in the setting of $C[a,b]$.)
Based on this inequality, Priestley \cite{priestley} has shown  that when $\{T_k\}$ is a sequence of positive linear transformations on ${\cal A}$ with $T_k(I)\leq I$ for all $k$, the set
$${\cal J}: = \{X\in {\cal A}: X^*=X,\,T_k(X)\rightarrow X,\,T_k(X^2)\rightarrow X^2\}$$
is a  norm-closed Jordan algebra of self-adjoint elements of ${\cal A}$, that is, a real linear
subspace of ${\cal A}$ closed under the Jordan product $X\circ Y:=\frac{XY+YX}{2}.$
We specialize Priestley's result by letting ${\cal A}={\cal M}_n$ (the space of all  $n\times n$ matrices) with $T_k=T$ for all $k$, where $T:{\cal M}_n\rightarrow {\cal M}_n$ is a positive unital linear transformation. Then, $\Hn$ is the Jordan algebra of self-adjoint elements of ${\cal M}_n$ (under the Jordan product mentioned above).
By the above result, the set ${\cal J}:=\{X\in \Hn: T(X)=X,\,\,T(X^2)=X^2\}$ is a Jordan subalgebra of $\Hn$. 
Suppose there is a $P\in \Hn$ with distinct eigenvalues such that $T(P)=P$ and $T(P^2)=P^2$. In $\Hn$, let  ${\cal E}$ be a Jordan frame with respect to which
$P$ has its spectral decomposition. Because $P$ has distinct eigenvalues, by uniqueness of Jordan frame (see \cite{faraut-koranyi}, Theorem III.1.1), this Jordan frame must be the Jordan frame of $P$ in the Jordan 
subalgebra ${\cal J}$. This means that ${\cal E}\subset {\cal J}$, showing that $T$ coincides with the Identity transformation on the span of ${\cal E}$. 
So, in the setting of $\Hn$, the implication $(a)\Rightarrow (b)$ in Theorem \ref{eja korovkin} can be deduced from Priestley's result. We note, however, that Priestley's result does not provide any information about $T^*$. 
Motivated by the above discussion,  we  raise the following  questions:
\begin{itemize}
{\it 
\item [$(1)$] Is there a  Kadison-type inequality for Euclidean Jordan algebras? That is, if $T$ is a positive linear transformation that is (sub)unital, can we assert that $T(x)^2\leq T(x^2)$ for all $x$? \\

%Note: In a private communication Jeong \cite{jeong}  has shown that $T(x)^2\underset{w}{\prec} T(x^2)$. 

\item [$(2)$] Does Priestley's result have an analog in the setting of Euclidean Jordan algebras? That is, if $T$ is a positive (sub)unital transformation on $\V$, can we say that the set
$$\{x\in \V: T(x)=x,\,\,T(x^2)=x^2\}$$
is a subalgebra of $\V$?
}
\end{itemize}

\gap

We now describe some examples where condition $(b)$ of Theorem  \ref{eja korovkin} holds.

\gap

\noindent{\bf Example 2} Let  $A=[a_{ij}]$ be an $n\times n$ real symmetric positive semidefinite matrix with every diagonal entry $1$ (that is, $A$ is a correlation matrix) and $\{e_1,e_2,\ldots, e_n\}$ be a Jordan frame in $\V$. Then, writing the Peirce decomposition of any $x\in \V$ as $x=\sum_{i\leq j}x_{ij}$, we define the  transformation
$$T: x\mapsto A\bullet x:=\sum_{i\leq j}a_{ij}x_{ij}.$$
This transformation is positive, unital, and self-adjoint (see Example 8 in \cite{gowda-doubly stochastic}) and satisfies condition $(b)$.
$T$ being doubly stochastic leads to some interesting consequences, see \cite{gowda-schur}. For example, by taking nonzero  numbers $a_1,a_2,\ldots, a_n$ and letting $A=\big [\frac{2a_ia_j}{a_i^2+a_j^2}\big ]$, we get the pointwise majorization inequality  
$[a_ia_j]\bullet x\prec \big [\frac{a_i^2+a_j^2}{2}\big ]\bullet x$. Written in the familiar form, this becomes
$$P_a(x)\prec L_{a^2}(x)\quad (x\in \V),$$
where for any $a=a_1e_1+a_2e_2+\cdots+a_ne_n\in \V$, $L_a(x):=a\circ x$ and $P_a(x):=2a\circ (a\circ x)-a^2\circ x$. 

\gap

\noindent{\bf Note:} In the case of $\Hn$, the real matrix $A$ can be modified 
as follows. Let $B$ be an $n\times n$ complex (Hermitian) positive semidefinite matrix with every diagonal entry $1$. Writing $X=[x_{ij}]$ and $B=[b_{ij}]$, we define the (Schur/Hadamard product) 
transformation $T$ on $\Hn$ by $T(X):=B\bullet X:=[b_{ij}x_{ij}]$. Then, $T$  is positive (by Schur product theorem, see \cite{horn-johnson}, Theorem 5.2.1) and
$T(E_i)=E_i$ for all $i$, where $E_i\in \Hn$ is the matrix with $1$ in $(i,i)$ slot and zeros elsewhere.
As $T^*(X)=\overline{B}\bullet X$, where $\overline{B}$ is the matrix of conjugates of entries of $B$, we see that $T^*(E_i)=E_i$ for all $i$. Note that, generally, $T$ need not be self-adjoint.
 
\gap

\noindent{\bf Example 3} Let $\V=\Hn$. We consider a {\it completely positive} linear transformation $T$ on $\Hn$ which, by definition, is of the form 
$$T(X):=A_1XA_1^*+A_2XA_2^*+\cdots+A_NXA_N^*\quad (X\in \Hn ),$$
where $A_1,A_2,\ldots, A_N$ are $n\times n$  complex matrices. 
If this transformation satisfies condition $(b)$ of the above theorem, then 
$T(I)=I=T^*(I)$ and so 
$A_1A_1^*+A_2A_2^*+\cdots+A_NA_N^*=I$ and
$A_1^*A_1+A_2^*A_2+\cdots+A_N^*A_N=I.$
We now characterize completely positive transformations satisfying condition $(b)$. \\
Let $T$ be as above and let 
$\{e_1,e_2,\ldots, e_n\}$ be a Jordan frame in $\Hn$ with $T(e_i)=e_i$ for all $i$. Let $\{E_1,E_2,\ldots, E_n\}$ denote the canonical Jordan frame in $\Hn$, where $E_i$ is a diagonal matrix with $1$ in the $(i,i)$ slot and zeros elsewhere. As $\Hn$ is simple, the Jordan frame $\{E_1,E_2,\ldots, E_n\}$ can be mapped into the Jordan frame $\{e_1,e_2,\ldots, e_n\}$ by an automorphism. Hence,  there is a unitary matrix $U$ such that 
$e_i=UE_iU^*$ for all $i$. Define the transformation $S$ on $\Hn$ by 
$$S(X):=U^*T(UXU^*)U\quad (X\in \Hn).$$
Then, $S$ is positive and $S(E_i)=E_i$ for all $i$. We see that 
$S(X)=B_1XB_1^*+B_2XB_2^*+\cdots+B_NXB_N^*,$ where $B_k=U^*A_kU$ for all $k$. 
By considering the block form of each $E_i$, we deduce from $S(E_i)=E_i$ that $B_k$ is  a diagonal matrix. Let $b_k$ denote the diagonal of $B_k$ (viewed as a column vector). Then, $B_kXB_k^*$ can be written as $(b_kb_k^*)\bullet X$. Letting $C:=\sum_{k=1}^{N} b_kb_k^*$, we see that
$$S(X)=C\bullet X\quad (X\in \Hn).$$
We observe that $C$, being a sum of rank-one matrices, is  positive semidefinite; from $S(E_i)=E_i$ we see that each diagonal entry of $C$ is one. Finally, 
\begin{equation} \label{schur form of cp transformation}
T(X)=U\Big(C\bullet U^*XU\Big )U^*\quad (X\in \Hn).
\end{equation}
Clearly, the above arguments can be reversed to see that a transformation of the  form (\ref{schur form of cp transformation}) is completely positive and satisfies condition $(b)$ of the theorem. \\
We now specialize by letting $k=1$. Let $T$ (defined by $T(X)=A_1XA_1^*$) coincide with the Identity transformation on the Jordan frame 
$\{e_1,e_2,\ldots, e_n\}$ in $\Hn$. Then, by the above, $B_1$ is a diagonal matrix where every diagonal entry has absolute value $1$ and 
 $A_1$ is  unitarily similar to $B_1$. 
%%%%%%%%%%%%%%%%%%%%%%%%%%%%%%%%%%%%%%%
\section{Sequential and weak-majorization versions}

Our next result deals with the sequential version of Theorem \ref{eja korovkin}. First, some preliminary material.
On the Euclidean Jordan algebra $\V$, for any $x\in \V$ with eigenvalues $x_1,x_2,\ldots,x_n$,  the $\infty$-norm is defined by 
$$||x||_\infty=\max_{1\leq k\leq n} |x_k|.$$
(It is known that $||\cdot||_\infty$ is a  norm on $\V$,  see e.g., \cite{gowda-doubly stochastic}.)
For any linear transformation $S$ on $\V$, let $||S||_\infty$ denote the operator norm relative to $||\cdot||_\infty$.
Now assume that $S$ is positive. For any $x\geq 0$, we have $0\leq x\leq ||x||_\infty\,e$ and so
$$0\leq S(x)\leq ||x||_\infty\,S(e).$$ We now apply (\ref{hirzebruch}) to get the inequality 
$||S(x)||_\infty\leq ||x||_\infty\,||S(e)||_\infty$ for all $x\geq 0$.
\\
Now let $x\in \V$. By considering the spectral decomposition of $x$, we can write $x=a-b$, where
$a=x^+$ and $b=x^-$ (see Section 2 for definitions).
Then, $||a||_\infty\leq ||x||_\infty$ and $||b||_\infty\leq ||x||_\infty.$ Hence,
$$||S(x)||_\infty=||S(a)-S(b)||_\infty\leq ||S(a)||_\infty+||S(b)||_\infty\leq ||a||_\infty\,||S(e)||_\infty+||b||_\infty\,||S(e)||_\infty$$
and so
$$||S(x)||_\infty\leq 2||x||_\infty\,||S(e)||_\infty.$$
Hence, for a positive linear transformation $S$ on $\V$,
$$||S||_\infty \leq 2||S(e)||_\infty.
$$
Now, let 
\begin{center}
${\cal B}(\V,\V)$:\,=\,Space of all  linear transformations from $\V$ to $\V$.
\end{center}
Since $\V$ is finite dimensional, the norm induced by the given inner product on $\V$ is equivalent to the $\infty$-norm. Correspondingly, the 
operator norms induced by these on ${\cal B}(\V,\V)$ are also equivalent. Hence, there is a positive constant $C$
(depending only on the dimension of $\V$) such that  for any positive linear transformation $S$ on $\V$,
\begin{equation} \label{norm of a positive transformation}
||S||\leq C||S(e)||,
\end{equation}
where $||S||$ is the operator norm of $S$ and $||S(e)||$ is the norm of $S(e)$ relative to the  norm induced by the given inner product on $\V$.

\begin{theorem}$\,\,$
{\it Suppose $\V$ is a Euclidean Jordan algebra of rank $n$ and $\{T_k\}$ is a sequence of positive linear transformations on $\V$ such that
$T_k(h)\rightarrow h$
for all $h\in \{e,p,p^2\}$, where $p\in \V$ with  distinct eigenvalues. Let
$p=p_1e_1+p_2e_2+\cdots+p_ne_n$ be the spectral decomposition of $p$. Then,
$$T_k(h)\rightarrow h\quad\mbox{and}\quad T_k^*(h)\rightarrow h$$
for all $h\in \mbox{span}\,\{e_1,e_2,\ldots, e_n\}$.
}
\end{theorem}

\gap

\begin{proof} 
For any $k$, $T_k$ is positive; hence, from (\ref{norm of a positive transformation}), 
$$||T_k||\leq C||T_k(e)||.$$
As $T_k(e)\rightarrow e$,  the sequence $||T_k(e)||$ is bounded. Hence, from the above,
the sequence $\{T_k\}$ is bounded in ${\cal B}(\V,\V)$.\\
We now  claim that $T_k(h)\rightarrow h$
for all $h\in \mbox{span}\,\{e_1,e_2,\ldots, e_n\}$. Since $T_k$s are linear, it is enough to show that $T_k(e_i)\rightarrow e_i$ for all $i=1,2,\ldots, n$. Suppose this is false;
assume, without loss of generality, that $T_k(e_1)\not\rightarrow e_1$. Then there is a subsequence $\{T_{k_l}\}$ of $\{T_k\}$ and a positive number $\varepsilon$ such that
\begin{equation}\label{greater than epsilon}
||T_{k_l}(e_1)-e_1||\geq \varepsilon\quad \mbox{for all}\,\,l.
\end{equation}
On the other hand, $T_{k_l}$ is a bounded sequence (in the finite dimensional space ${\cal B}(\V,\V)$), hence has a subsequence -- continue to call this $T_{k_l}$ -- that converges to a linear transformation, say, $T$.
We see that $T$ is positive and (by the imposed conditions on $\{T_k\}$) satisfies the conditions
$$T(e)=e,T(p)=p,\,\,\mbox{and}\,\,T(p^2)=p^2.$$
Now, by Theorem \ref{eja korovkin}, $T(e_i)=e_i$ for all $i=1,2,\ldots, n$. In particular, $T(e_1)=e_1$. But this means that $T_{k_l}(e_1)\rightarrow e_1$ contradicting (\ref{greater than epsilon}).
Hence, $T_k(e_i)\rightarrow e_i$ for all $i=1,2,\ldots, n$.
\\
We now claim that $T_k^*(e_i)\rightarrow e_i$ for all $i$.  
Suppose, without loss of generality, $T_k^*(e_1)\not\rightarrow e_1$. Since $||T_k^*||=||T_k||$ for all $k$,  the sequence $\{T_k^*\}$ is  bounded in ${\cal B}(\V,\V)$. 
Then, as argued before, there is a subsequence $T_{k_l}^*$ such that 
$T_{k_l}^*(e_1)\not\rightarrow e_1$ and $T_{k_l}^*$ converges to, say, $S$. As the adjoint operation is continuous, we have $T_{k_l}\rightarrow S^*$. As  $T_{k_l}(e_i)\rightarrow e_i$ for all $i$, we have $S^*(e_i)=e_i$ for all $i$. Since $S^*$ is positive, from our previous result, $S(e_i)=e_i$ for all $i$. But then, $T_{k_l}^*(e_i)\rightarrow e_i$ for all $i$, contradicting our  assumption that $T_{k_l}^*(e_1)\not\rightarrow e_1$. Thus, we have our claim.   

\end{proof}

\gap

We now state a result that is analogous to Theorem \ref{weak majorization korovkin} on a  {\it simple} algebra. Recall that in $\V$, by definition, $u\underset{w}{\prec}v$ if 
$\lambda(u)\underset{w}{\prec}\lambda(v)$  in $\Rn$.

\gap

\begin{theorem} \label{eja weak majorization korovkin}$\,\,$
{\it
Let $\V$ be a simple Euclidean Jordan algebra of rank $n$ and $T:\V\rightarrow \V$  be a positive linear transformation. 
Let $p\in \V$ with spectral decomposition $p=p_1e_1+p_2e_2+\cdots+p_ne_n$, where $p_i$s are  positive and distinct. 
Suppose the following conditions hold:
\begin{itemize}
\item [(i)] $T(e)\underset{w}{\prec} e$, $p\underset{w}{\prec} T(p)$, and $T(p^2)\underset{w}{\prec}p^2$, and 
\item [(ii)] $T(p)$ and $T(p^2)$ operator commute.
\end{itemize}
Then 
\begin{itemize}
\item [(a)] $T$ coincides with an automorphism on $\mbox{span}\{e_1,e_2,\ldots, e_n\}$, 
\item [(b)] $\lambda(T(x))=\lambda(x)\,\,\mbox{for all}\,\,x\in \mbox{span}\{e_1,e_2,\ldots, e_n\},$ and 
\item [(c)]$T$ is doubly stochastic.
\end{itemize}
}
\end{theorem}

\gap

\begin{proof}
We assume that all the assumptions are in place. 
By permuting $e_1,e_2,\ldots, e_n$, we may assume that 
$p_1>p_2>\cdots>p_n$. Then, $p=\lambda(p)*{\cal E}$, where ${\cal E}:=\{e_1,e_2,\ldots,e_n\}$. Since $T(p)$ and $T(p^2)$ operator commute, they 
have their  spectral representations with respect to the same Jordan frame, say, ${\cal F}=\{f_1,f_2,\ldots, f_n\}$. We write $T(p)=r*{\cal F}=r_1f_2+r_2f_2+\cdots+r_nf_n$ and $T(p^2)=s*{\cal F}:=s_1f_1+s_2f_2+\cdots+s_nf_n$, where $r=(r_1,r_2,\ldots,r_n)$ and $s=(s_1,s_2,\ldots, s_n)$; we assume, without loss of generality, that the entries of $r$ are decreasing. Note that 
$r=r^\downarrow=\lambda(T(p))$ and $s^\downarrow=\lambda(T(p^2))$. 
Now, since  $\V$ is simple, there is an  automorphism $\phi$ which takes ${\cal F}$ to ${\cal E}$, so $\phi(f_i)=e_i$ for all $i$, see \cite{faraut-koranyi}, Theorem IV.2.5. Then,
$$\phi(T(p))=\phi(r*{\cal F})=r*{\cal E}\,\,\mbox{and}\,\,\phi(T(p^2))=\phi(s*{\cal F})=s*{\cal E}.$$ 
Let $S:=\phi\circ T$ and $\overline{p}:=\lambda(p)$. Then, $S$ is positive and 
$$T(e)\underset{w}{\prec} e\Rightarrow \lambda(T(e))\underset{w}{\prec}\lambda(e)\Rightarrow \lambda(T(e))\leq \lambda(e)\Rightarrow T(e)\leq e\Rightarrow S(e)\leq e,$$
where the second implication is due to the fact that $\lambda(e)$ is the vector of $1$s in $\Rn$. Now consider 
the  matrix $B$ defined by
$$B=[b_{ij}],\quad b_{ij}:=\frac{1}{||e_i||^2}\langle S(e_j),e_i\rangle.$$
Since $S$ is positive and $S(e)\leq e$, we see that 
$B$ is nonnegative  and $\sum_{j=1}^{n}b_{ij}\leq 1$ for all $i$. 
From the relations $$S\Big (\sum_{j=1}^{n}\bar{p}_je_j\Big )=S(p)=r*{\cal E}=\sum_{i=1}^{n}r_ie_i\quad\mbox{and}\quad
S\Big (\sum_{j=1}^{n}\bar{p}_j^2e_j\Big )=S(p^2)=s*{\cal E}=\sum_{i=1}^{n}s_ie_i$$
we verify that $B\bar{p}=r$ and $B\bar{p}^2=s.$
Moreover,  from condition $(i)$, as $\phi$  preserves eigenvalues,
$$\lambda(p)\underset{w}{\prec}\lambda(T(p)) =r=\lambda(S(p))\,\, \mbox{and}\,\, \lambda(S(p^2))=s^\downarrow\underset{w}{\prec}\lambda(p^2).$$
In summary:  
$B$ is nonnegative, subunital, $\overline{p}\underset{w}{\prec} B\overline{p}$, and 
$B\overline{p}^2\underset{w}{\prec}\overline{p}^2;$ additionally, the entries of $\overline{p}$ are strictly decreasing and those of $B\overline{p}$ are decreasing.  
\\
From  Theorem \ref{weak majorization korovkin}, we see that $B$ is the Identity matrix. So, for all $i,j$,
$$\langle S(e_j),e_i\rangle =||e_i||^2\delta_{ij}.$$
From Lemma \ref{basic lemma}, $S=I$ on $\W:=\mbox{span}\{e_1,e_2,\ldots, e_n\}$. 
So,
$\phi(T(x))=x$ for all $x\in \W$, that is,
$$T(x)=\phi^{-1}(x) \,\, \mbox{for all}\,\,x\in \W.$$
As $\phi^{-1}$ is an automorphism on $\V$, we  have Item $(a)$.
Since automorphisms preserve eigenvalues, $T$ preserves eigenvalues of every element in $\W$. This gives $(b)$.
Finally, $(c)$ comes from Corollary \ref{T equals  an  automorphism}.
\end{proof}

\gap

\noindent{\bf Remarks.} We note that conditions $(i)$ and $(ii)$ in the above result are necessary and sufficient for $T$ to coincide with an automorphism 
on $\{e_1,e_2,\ldots, e_n\}$. Moreover, in the presence of $(ii)$, $(i)$ is equivalent to each of the following:
\begin{itemize}
\item [$(1)$] $\lambda(T(h))=\lambda(h)$ for all $h\in \{e,p,p^2\}$.
\item [$(2)$] $T(e)\underset{w}{\prec} e$, $\lambda(T(p))=\lambda(p)$,  
and $T(p^2)\underset{w}{\prec}p^2$. 
\end{itemize}
It is not clear if the assumption that $\V$ is simple can be dispensed with. 

\gap

\noindent{\bf Acknowledgments and concluding remarks}
Thanks are due to Michael Orlitzky, Roman Sznajder, and Juyoung Jeong for their comments and suggestions.
In a private communication \cite{jeong-private}, Jeong notes that Theorem 3.1 and Theorem 4.2 continue to hold when the quadratic function $t\mapsto t^2$ on $\R$ is replaced by a strictly convex function. He also shows (by an example) that when $n\geq 4$, the set $\Omega_p$ (that appears in the problem posed in Section 3) may  contain matrices other than doubly stochastic ones. 

%%%%%%%%%%%%%%%%%%%%%%%%%%%%%%%%%%%%%%%%%%%%%%%%%%%%%%%%%%%%%%%%%%%%%%%%%%%%%%%%

\end{document}